\documentclass{article} 

\usepackage{mathmag}

\theoremstyle{theorem}

\theoremstyle{definition}

\allowdisplaybreaks

\makeatletter
\@addtoreset{footnote}{page}
\makeatother

\begin{document}

\title{The Irrationality of $e$ and $\pi$}
	
\author{L. Lerner\\     %%%%  DO NOT FILL IN AUTHOR'S NAMES
\scriptsize School of Physical, Chemical and Earth Sciences\\           %%%%  UNTIL YOU RECIEVE YOUR
Adelaide University, Adelaide, Australia\\                       %%%%  PROVISIONAL ACCEPT LETTER.
leonid.lerner@adelaide.edu.au}                            
	
\maketitle

%\begin{abstract}

%\end{abstract}

\section{Introduction}

Unlike an ordinary fraction with an infinite numerator and denominator, an infinite continued fraction must be irrational \cite{stein}. Euler \cite{euler} was the first to show that the base of the natural logarithm $e$ is irrational, by numerically estimating its continued fraction, and showing the infinite sequence of convergents $p_n/q_n$ so obtained converged to $e$ using the Ricatti equation. Hermite \cite{hermite} showed that the recurrence relations for these convergents correspond to the recurrence relations between certain improper integrals, so proving the continued fraction tends to $e$ in the limit of infinite $n$. Here we provide a motivation for the integrals involved, complimentary to that in \cite{cohn}, and obtain closed form integral representations for $p_n$ and $q_n$ for all $n$. An interesting feature, is that the above results provide a short cut to the standard proof of the irrationality of $\pi$.

\section{The Irrationality of $e$}
Carrying out Euclid's algorithm gives
\begin{eqnarray*}
e =
2+\cfrac{1}{1+\cfrac{1}{2+\cfrac{1}{1+\cfrac{1}{1+
				\cfrac{1}{4+
							\cdots\vphantom{\cfrac{1}{1}} }}}}}
\equiv
 [2, 1, 2, 1, 1, 4, 1, 1, 6, 1, 1, 8, \ldots]. 
\end{eqnarray*}
To prove that $e$ is irrational we have to show that the above pattern, obtained numerically, actually follows from the definition of $e$. Writing the $n$\textsuperscript{th} convergent as $p_n/q_n=[a_0,a_1,\ldots,a_n]$ where $p_n,q_n$ are integers, one can show by induction\footnote{We note the replacement $a_n \rightarrow a_n + 1/a_{n+1}$ relates two successive convergents. Then by induction
\begin{eqnarray*}
	\frac{p_{n+1}}{q_{n+1}} = \frac{p_{n-2} + (a_n+1/a_{n+1}) \cdot p_{n-1}}{q_{n-2} + (a_n+1/a_{n+1}) \cdot q_{n-1}}	
	= \frac{p_{n} + p_{n-1}/a_{n+1}}{q_{n} + q_{n-1}/a_{n+1}} 
	= \frac{p_{n-1} + a_{n+1} \cdot p_{n}}{q_{n-1} + a_{n+1} \cdot q_{n}}
\end{eqnarray*}} 
that 
\begin{eqnarray*}
p_n = p_{n-2} + a_n \cdot p_{n-1}, \qquad
q_n = q_{n-2} + a_n \cdot q_{n-1}
\end{eqnarray*}
where the starting point is $(p_{-1},  p_0) = (1, a_0)$ and $(q_{-1},  q_0) =  (0, 1)$.

We now specialize to $\mathrm{e}= [2, 1, 2, 1, 1, 4, 1, 1, 6, \ldots]$. 
With the exception of the first two elements, the sequence $a_n$ has a repeat pattern of length 3 in the sense that:
$$ (a_{3i-1}, a_{3i}, a_{3i+1} )= (2  i, 1, 1) $$
This is not a numerically repeating sequence, for then $e$ would be the solution to a quadratic, but rather an algebraic sequence.
Because the equation for $a_n$ depends on $n\! \mod 3$, the equations for $p_n$ and $q_n$ also have this dependence. We can remove this by reordering $p_n$ and $q_n$ into 3 new sequences $A_i$, $B_i$, $C_i$ as
\begin{eqnarray*}
(p_{3i-1}, p_{3i}, p_{3i+1} ) &\equiv& (B_i, C_i, A_i) \\
(q_{3i-1}, q_{3i}, q_{3i+1} ) &\equiv& (B_i', C_i', A_i')
\end{eqnarray*}
We find the linear relations between the sequences $A_i$, $B_i$, $C_i$ by substituting the above definitions into the recursion relation for $p_n$ and $q_n$. Thus
\begin{eqnarray*}
B_i &=& 2  i  A_i    + C_{i-1} \\
A_i &=&            C_{i-1} + B_{i-1} \\
C_i &=&            A_i     + B_i
\end{eqnarray*}
We can eliminate $B_i$, and $C_i$ to find the relationship between the $A_i$'s
$$ A_{n+1} = 2  (2  n + 1) A_n + A_{n-1}. $$
The sequences $p_n$ and $q_n$, while obeying this recurrence relation, differ only in the initial conditions. It is most convenient to use $A_0$ and $A_1$ for this, because they are smallest. For the sequence $p_n$ these are $(A_0, A_1) = (1, 3)$, and for the sequence $q_n$ they are $(A_0', A_1') = (1, 1)$. Although $p_{-2}$ which corresponds to $A_0$ does not appear in a convergent of $e$, it can be deduced from the relation $A_2 = 6  A_1 + A_0$, by substituting for $A_1$ and $A_2$, which correspond to $p_1$ and $p_4$ respectively. And similarly for $q_{-2}$.

We now try solving the equation for $A_n$ by finding an integral with the same recurrence relation. For instance, integrating $ a_n = \int_0^\infty  u^n e^{-u} du $ by parts
gives $a_n = n  a_{n-1}$ with $a_0=1$. Solving, we get get $a_n = n!$. In this case reverting the integrand to the same form in a single partial integration, was responsible for a recurrence relation of first order.
To get a second-order relation we need to find an integrand which is reverted to the same form not by one, but by two integrations by parts. Replacing $u^n$ by an arbitrary function $f(u)$ shows that $f''(u)$ should have the same functional form as $f(u)$, while $f'(u)$ should not. Choosing $f(u)=(1-u^2)^n$ we arrive at the integrand 
$$t_n(u) = (1 - u^2)^n e^{-u},$$
which satisfies our criterion. A single integration by parts changes this integrand to $-2 n u (1-u^2)^{n-1}  e^{-u}$, which is not of the original form, while a further partial integration gives 
$$-2 n  (1-u^2)^{n-1}\ e^{-u} +4  n  (n-1)  u^2
(1-u^2)^{n-2} e^{-u}. $$
The first term is just $- 2nt_{n-1}$, while the identity $u^2 = 1 - (1 - u^2)$, changes the second term to $4  n (n-1)  (t_{n-2}-t_{n-1})$. The complete integrand is now expressible in terms of $t_n(u)$, and we get a recurrence relation. The range of integration needs to be such that the boundary term in the partial integration is zero, as is the case for $x^n e^{-x}$ on the range $[0, \infty]$. This ensures no additive constant is introduced by the boundary terms on each partial integration, so that the recurrence relations are homogeneous. The integrand $(1-u^2)e^{-u}$ vanishes at the points $u = {-1, 1, \infty}$. Choosing $[1, \infty]$  as our integration range, we define
$$ T_n = \int_1^\infty (1 - u^2)^n e^{-u/a} du$$
which has the recursion relation
$$T_{n+1} 	= -2 a^2 (n+1) (2n + 1) T_n + 4 a^2 n(n+1) T_{n-1}.$$
This differs from the recurrence relation for $A_n$ only in the constants on the right hand side. These can be brought into the correct form by transforming to a new sequence related to the original by $S_n = T_n/(b^n n!)$, where $b$ is a constant to be chosen. Substituting for $T_n$ we get a recurrence relation for the $S_n$
$$S_{n+1} (n+1)! b^{n+1} = -2 a^2 n! b^n (n+1) (2n + 1) S_n + 4 a^2 b^{n-1} n(n+1) S_{n-1} $$
Dividing both sides by $(n+1)! b^{n+1}$, and setting $a = 2$ and $b = -4$ we get
$$S_{n+1} = 2 (2 n + 1) S_n + S_{n-1}$$
which is the recurrence relation for $A_n$ with 
$$ S_n = \frac{1}{ 4^n n!} \int_1^\infty \textrm{d}u (u^2 - 1)^n e^{-u/2}. $$
The initial values for the recurrence relation are 
$S_0 = 2 e^{-1/2}$, 
$S_1 = 6 e^{-1/2}$.
Comparing to $A_0 = 1, A_1 = 3$ we see the ratio of the two initial values of the sequences $A_n$ and $S_n$ are the same, so that $A_n$ is a multiple of $S_n$. Thus
$$ A_n = 2^{-1}e^{1/2} S_n, $$
and we have an integral representation for $p_{3n+1}$
$$ p_{3n+1} = \frac { e^{1/2}} {2^{2n+1} n!} \int_1^\infty \textrm{d}u (u^2 - 1)^n e^{-u/2}.$$
For $q_{3n+1}$ we need to find an independent integral obeying the same recurrence relations with the initial values $A_0' = 1, A_1' = 1$. We can then use a linear combination with $A_i$ to satisfy the initial values of $q_n$. This suggest the same integrand, but over $u = [-1, \infty]$, where the integrand also vanishes at both integration limits. Defining
$$ S_n' = \frac { e^{1/2}} {2^{2n+1} n!} \int_{-1}^\infty \textrm{d}u (u^2 - 1)^n e^{-u/2}$$
we have
$S_0' = S_1' =e$.
We now see that $S_1'/S_0' = A_1'/A_0'$ and therefore $S_n'$ is a multiple of $A_n'$. In particular $S_n' = e A_n'$, and we have an integral representation of $q_{3n+1}$
$$
	q_{3n+1} 
	= \frac { e^{-1/2}} {2^{2n+1} n!} \int_{-1}^\infty \textrm{d}u (u^2 - 1)^n e^{-u/2}.   
$$   
and so the convergent of our continued fraction becomes 
$$ p_{3n+1}/q_{3n+1} = \left. e \int_{1}^\infty \textrm{d}u (u^2 - 1)^n e^{-u/2} \middle / \int_{-1}^\infty \textrm{d}u (u^2 - 1)^n e^{-u/2} \right. . $$
Now since both $p_{3n+1}$ and $q_{3n+1}$ tend to infinity with $n$, it follows that the integrals in the numerator and the denominator above also both tend to infinity. Their difference however is bounded, 
$$
 \int_{-1}^1 \textrm{d}u (1-u^2 )^n e^{-u/2}  < 2 \, e^{1/2}
$$
which means their ratio tends to unity.
\footnote{One can get an asymptotic limit for the integral as $n \rightarrow \infty$ by expanding the exponential and changing variables $t=u^2$ which gives a sum of $\beta$-functions
$$\sum_{m=0}^\infty  4^{-m} (2m!)^{-1} \int_0^1 \textrm{d}t  t^{m-1/2}(1-t)^n =
\sum_{m=0}^\infty 4^{-m} (2m!)^{-1}B(n+1,m+1/2) \sim \sqrt{\pi/(n+1)}.
$$
}
We then have 
$$ \lim_{n \to \infty} \frac{p_{3n+1}}{q_{3n+1}} = e $$
Moreover this holds for all $p_n/q_n$ because the simple continued fraction is always convergent. Alternatively, from the relations between $A_n, B_n, C_n$, it follows that $\lim_{n \rightarrow \infty} B_n =  2nA_n$, so $\lim_{n \rightarrow \infty}B_n/B_n' = A_n/A_n' = e$. Then, from  $C_n = A_n + B_n$ it follows that $\lim_{n \rightarrow \infty} C_n/C_n'  \rightarrow e$. 

Thus the infinite continued fraction $[2, 1, 2, 1, 1, 4, 1, 1, 6, \ldots]$ equals $e$, from which it follows that $e$ is irrational. We can also evaluate the relative error $r_n=|p_{3n+1}/ q_{3n+1} - e |/e$ giving
\begin{eqnarray*}
r_n = \int_{-1}^1 \textrm{d}u (1 - u^2)^n e^{-u/2} /  \int_{-1}^\infty \textrm{d}u (u^2 - 1)^n e^{-u/2} 
 <  \frac{ e^{1/2}}{4^{n} (2n)!- e^{1/2}}.
\end{eqnarray*}
since
\begin{eqnarray*}
\int_{1}^\infty \textrm{d}u (u^2 - 1)^n e^{-u/2} 
<\int_0^\infty \textrm{d}u \, u^{2n} e^{-u/2} = 2^{2n+1}(2n)!.
\end{eqnarray*}

\section{The Irrationality of ${\pi}$} The integral
$$R_n= \Re \int_{-1}^1 \textrm{d}u (1-u^2 )^n e^{i xu} =\int_{-1}^1 \textrm{d}u (1-u^2 )^n \cos xu $$
has the same recurrence relation as $T_n$ with $a=1/ix$, since  the latter does not mix real and imaginary components
$$x^2 R_{n+1} 	= 2 (n+1) (2n + 1) R_n - 4 n(n+1) R_{n-1}.$$
We get a more useful form by substituting $R_n=  2^{n+1} x^{-(2n+1)}  n!  Q_n$, so the recurrence relation for the new sequence $Q_n$ reads
$$Q_{n+1} 	=  (2n + 1) Q_n -  x^2 Q_{n-1},$$
with the inital values $Q_0=\sin x, Q_1=\sin x -x \cos x$. Now set $x=\pi/2$. It follows immediately by inductance on the recurrence relation that $Q_n$ is a polynomial in $\pi/2$ with degree $\le n$ and with integer coefficients. Then if $\pi$ is rational, $\pi/2=a/b$ for some integers $a,b $ and we have 
$$ b^{2n+1}Q_n= \frac{a^{2n+1}}{2^{n+1} n!}\int_{-1}^1 \textrm{d}u (1-u^2 )^n \cos xu,$$
where the LHS is an integer on account of $d(Q_n) \le n$. Moreover the RHS is positive definite and so $b^{2n+1}Q_n\ge 1$. However the integral on the RHS is bounded, and so the RHS vanishes in the limit $n \to \infty$. This is inconsistent with the limit on the LHS, and so a rational $\pi$ leads to a contradiction.

%\begin{biog}
%\item[Author Name 1] Insert author bio here.
%\begin{affil}
%Department of Mathematics, University America, Washington DC 20036\\
%authorname@ua.edu
%\end{affil}
%\end{biog}

%\begin{biog}
%\item[Author Name 2] Insert author bio here.
%\begin{affil}
%Department of Mathematics, University America, Washington DC 20036\\
%authorname@ua.edu
%\end{affil}
%\end{biog}

\vfill\eject

\end{document}